# HOW A STRONGLY IRREDUCIBLE HEEGAARD SPLITTING INTERSECTS A HANDLEBODY

MATT JONES AND MARTIN SCHARLEMANN

ABSTRACT. In [Sc1] it is shown that a strongly irreducible Heegaard splitting surface $Q$ of a 3-manifold $M$ can, under reasonable side conditions, intersect a ball or a solid torus in $M$ in only a few possible ways. Here we extend those results to describe how $Q$ can intersect a handlebody in $M$.

## 1. INTRODUCTION AND PRELIMINARIES

A central problem in the century-old theory of Heegaard splittings is to understand how different Heegaard splittings of the same 3-manifold might compare. Much progress on this and other questions has followed the introduction, by Casson and Gordon, of the idea of *strongly irreducible* Heegaard splittings (see [CG], [Sc2]). For example, in [RS1] it was shown that two strongly irreducible Heegaard splittings of the same closed orientable 3-manifold can be isotoped to intersect each other in a particularly helpful way. A more general question, useful in the same endeavor, is to understand how a strongly irreducible splitting surface intersects a handlebody lying in the interior of the 3-manifold. That is the goal of this paper.

In general the intersection of the splitting surface and a handlebody can be quite complicated, even when the splitting is strongly irreducible. But, with reasonable side-conditions on how the splitting surface intersects the handlebody boundary, the picture has been shown to simplify dramatically when the handlebody is either a ball or a torus [Sc1]. Here we obtain similar results (but not as precise) when the handlebody is of arbitrary genus. The core of our approach is to examine just the part of the splitting surface that lies within the handlebody. It turns out that, under reasonable side conditions, this surface is *weakly incompressible* (see Definition 1.1). So the general problem translates directly into the problem of classifying weakly incompressible surfaces embedded in a handlebody.

Although the results here are fairly technical, they may also prove useful. For example, Theorem 2.1 is helpful in classifying how two genus two Heegaard splittings of the same manifold intersect (cf. [RS3]).

All 3-manifolds will be orientable. A compression body $H$ is a 3-manifold obtained from $(surface) \times I$ by attaching 2-handles to $(surface) \times \{1\}$ and capping off any

*Date*: February 1, 2018.

1991 *Mathematics Subject Classification.* 57N10.

First author partially supported by a grant from the UCSB College of Creative Studies. Second author partially supported by an NSF grant.

2-sphere boundary components that result with 3-balls. $\partial_+ H$ is $(surface) \times \{0\}$ and $\partial_- H = \partial H - \partial_+ H$. (Mnemonic: $\partial_+ H$ is a more complicated surface then $\partial_- H$.)

The compression body $H$ is a *handlebody* if $H$ is a compression body with $\partial_- H = 0$. A *splitting suface* for a *Heegaard splitting* of a 3-manifold $M$ is a properly embedded closed orientable surface $S$ which divides $M$ into compression bodies $H_1$ and $H_2$ so that $M = H_1 \cup_S H_2$ and $\partial_+ H_1 = S = \partial_+ H_2$. The splitting is *weakly reducible* if there is a disjoint pair of compression disks for $S$, one in $H_1$ and one in $H_2$. Otherwise, it is *strongly irreducible*. Here is a generalization of those ideas.

**Definition 1.1.** *A properly imbedded oriented surface $(Q, \partial Q) \subset (M, \partial M)$ is a splitting surface if $M$ is the union of two 3-manifolds $X$ and $Y$ along $Q$ so that $\partial X$ induces the given orientation on $Q$ and $\partial Y$ induces the opposite orientation. A compressing disk for $Q$ in $X$ (resp. $Y$) is called a* meridian disk *in $X$ (resp. $Y$) and its boundary a meridian curve for $X$ (resp. $Y$).*

*The splitting surface $Q$ is* bicompressible *if $Q$ is compressible into both $X$ and $Y$. $Q$ is called* strongly compressible *if there are meridian disks in $X$ and $Y$ with disjoint boundaries. If a splitting surface is not strongly compressible then it is called* weakly incompressible.

**Remarks:** If $Q$ is both bicompressible and weakly incompressible, there are compressing disks in both $X$ and $Y$, but any pair of such disks, one in $X$ and the other $Y$, necessarily have boundaries that intersect along $Q = X \cap Y$.

The following is a slight adaptation of a lemma and proof originally given in [Sc1].

**Lemma 1.2.** *Let $Q$ be a bicompressible, weakly incompressible properly imbedded surface in a handlebody $H$, dividing $H$ so that $H = X \cup_Q Y$. Let $F$ be obtained from $Q$ by a series of compressions into $X$, and let $\Gamma$ be the set of arcs dual to these compressions. Let $D$ be a disk or collection of disks, in general position with respect to $\partial Q$ and and $Q$ and with $\partial D \subset (\partial H \cup F) - \eta(\Gamma)$. Assume that $int(D) \cap F$ contains no closed curves. Then $\Gamma$ can be made disjoint from $D$ by a series of edge slides and isotopies.*

**Proof:** Details of the proof, which we here sketch, can be found in [ST, Prop. 2.2]. Choose a representation of $\Gamma$ which minimizes $|D \cap \Gamma|$, and assume that $D \cap \Gamma \neq \emptyset$. Choose a compressing disk $E$ for $Q$ in $Y$ which minimizes $|D \cap E|$. Note that if there were any closed curves of $D \cap E$ which bounded disks in $D$ disjoint from $\Gamma$, we could choose one innermost in $D$ and replace the disk it bounded in $E$ with a copy of the one it bounded in $D$ to reduce $|D \cap E|$.

Call an arc of $D \cap E$ which has both endpoints on the neighborhood of the same point of $D \cap \Gamma$ a *loop*. If there were a loop of $D \cap E$ which bounded in $D$ a disk disjoint from $\Gamma$, we could choose an innermost such loop and $\partial$-compress $E$ with the disk bounded by that loop. At least one of the two resulting disks would have to be essential, since their sum is $E$, so choosing that disk instead of $E$ would have given a lower value for $|D \cap E|$ (the chosen loop is eliminated, if nothing else).

By considering a component of $D \cap E$ which is innermost in $D$ among all closed curves and loops (if there are any such components), and by considering the disk in $D$ bounded by that component, we see that there must be some point $w \in D \cap \Gamma$ which is incident to no loops. But since $\eta(w) \cap D$ is a compressing disk for $Q$ in $X$, it must intersect $E$. So choose an arc $\alpha$ which is outermost in $E$ among all arcs of $D \cap E$ which are incident to $\eta(w)$. Then $\alpha$ cuts off from $E$ a disk $E'$ with $E' - \alpha$ disjoint from $w$. Let $e$ be the edge of $\Gamma$ which contains $w$. Then the disk $E'$ gives instructions about how to isotope and slide the edge $e$ until it is disjoint from $D$, and these slides and isotopies ultimately decrease $|D \cap E|$. □

**Lemma 1.3.** *The only closed splitting surface that is bicompressible, weakly incompressible and lies in a ball is an unknotted torus.*

**Proof:** Suppose $Q \subset B$ is a counterexample, dividing a ball $B$ into two parts $X$ and $Y$, with $\partial B \subset X$. Since $Q$ is weakly incompressible, all compressing disks must lie on the same component $Q_0$. Let $Q_X, Q_Y$ be the surfaces (possibly spheres) obtained by maximally compressing $Q_0$ into $X$ and $Y$ respectively and $W \subset B$ be the 3-manifold, containing $Q_0$, that lies between $Q_X$ and $Q_Y$. Then $Q_0$ is a Heegaard splitting surface for $W$ and, since $Q$ is weakly incompressible, the splitting is strongly irreducible. It follows from [CG] that $\partial W = Q_X \cup Q_Y$ is incompressible in $W$. By construction $\partial W$ is also incompressible in the complement of $W$, since a further compression would either be disjoint from $Q - Q_0$ (which violates maximality in the construction) or it must intersect $Q - Q_0$, forcing $Q$ to be compressible in the complement of $Q_0$, and so violating weak incompressibility of $Q$. Hence $\partial W$ is a collection of spheres. It follows easily that $Q = Q_0$ and this surface is a Heegaard splitting surface of $B$. By Waldhausen's theorem (see [ST] for an updated proof) any such splitting is standard. But any standard splitting of genus greater than one is strongly compressible. □

The next lemma is a variation of one originally given by Frohman in [Fr, Lemma 1.1].

**Lemma 1.4.** *Let $Q$ be a bicompressible, weakly incompressible splitting surface in a handlebody $H$, dividing $H$ into $X$ and $Y$. Suppose there is an essential loop $\alpha$ in $Q$ which intersects a compressing disk $C$ for $Q$ exactly once. If $\alpha \cup C$ is contained in the interior of a ball properly embedded in $H$, then $C$ is the meridian and $\alpha$ the longitude of an unknotted torus component of $Q$ that lies in a ball in $H$.*

**Proof:** We may as well assume that $Q$ is connected. Suppose, for $\alpha$ and $C$ as described, that $\alpha \cup C$ is contained in a ball $B$ and that $X$ is the side that contains $C$. Among all possible choices for $B$ and choices of curve, disk pairs $\alpha, C$ in $Q, X$ that intersect in a single point, chose $B, \alpha, C$ to minimize $\partial B \cap Q$. We will show that in fact $\partial B \cap Q = \emptyset$ so that Lemma 1.3 applies as required.

Indeed, we will show that if $\partial B \cap Q \neq \emptyset$ then this forces a contradiction. Suppose first that there were a component of $\partial B \cap Q$ that is inessential in $Q$. Choose an innermost one in $Q$ and let $E$ be the disk in $Q$ that it bounds, so that either $E \subset B$ or $E \subset H - B$. Furthermore, $E$ must be disjoint from $\alpha \cup C$, because $\partial E \subset \partial B$ is disjoint from $\alpha \cup C$ and $\alpha$ and $\partial C$ are both essential in $Q$. If $E \subset H - B$, the disk in $\partial B$ which is bounded by $\partial E$ could be replaced by a copy of $E$ to reduce $|\partial B \cap Q|$ (remember that handlebodies are irreducible). But if $E \subset B$, then cutting $B$ on $E$ produces two balls, and the one which contains $\alpha \cup C$ would have a lower value for $|\partial B \cap Q|$, again a contradiction. We conclude from these contradictions that each component of $\partial B \cap Q$ must be essential in $Q$.

If $|\partial B \cap Q| = 1$ then the disk of $\partial B - Q$ which lies in $Y$ would be a compressing disk for $Q$ which is disjoint from $C$, contradicting weak incompressibility of $Q$. So we can assume $|\partial B \cap Q| > 1$. Then there is a component $\beta$ of $\partial B \cap Q$ which bounds a disk $D$ in $\partial B$ with $int(D) \cap Q$ a non-empty collection of loops which are all innermost in $\partial B$. Let $P$ be the non-disk planar component of $D - Q$. The disks of $D - P$ are all compressing disks for $Q$ and are disjoint from $C$, so they must be in $X$, and thus $P \subset Y$. Let $F$ be the surface obtained from $Q$ by simultaneous compression on all these disks, and let $\Gamma$ be the set of arcs dual to these compressions. Now use Lemma 1.2 to isotope $\Gamma$ off of $D$. The slides involved should be done without dragging along $\partial C$ or $\alpha$ (so, technically, $\partial C$ and $\alpha$ may change as curves in $Q$). But the result of the slides would again violate the condition that $\partial B \cap Q$ had been minimized. $\square$

**Corollary 1.5.** *Suppose $Q$ is a bicompressible, weakly incompressible splitting surface in a handlebody $H$ and $Q$ can be obtained from a surface $F$ by attaching tubes along a set of arcs $\Gamma$. If there is a circuit in $F \cup \Gamma$ which is contained in the interior of a ball and which passes over some arc of $\Gamma$ exactly once, then a component of $Q$ is an unknotted torus in a ball in $H$.*

Such a circuit will be referred to as a cycle in a ball.

## 2. Characterizing bicompressible, weakly incompresible splitting surfaces

**Theorem 2.1.** *Let $Q$ be a properly embedded bicompressible weakly incompressible splitting surface in a handlebody $H$, dividing $H$ into 3-manifolds $X$ and $Y$. Suppose no component of $Q$ is an unknotted torus in a ball in $H$. Let $\Delta$ be a complete set of meridian disks for $H$. Then there is a properly imbedded incompressible surface $F \subset H$ intersecting $\Delta$ only in arcs and a set $\Lambda$ of pairwise disjoint arcs so that $\Lambda \subset \Delta$, $\partial \Lambda = \Lambda \cap F$ and $Q$ is properly isotopic to the surface obtained from $F$ by attaching tubes along the arcs of $\Lambda$. The side, $X$ or $Y$, containing the interiors of these tubes can be specified in advance.*

**Proof:** According to Corollary 1.5, there is no cycle in a ball. Begin by maximally compressing $Q$ into $X$, say. At the end of this we are left with a surface $F$ which divides $H$ into $X_-$ and $W$, where $X_-$ is the boundary-reduced $X$ and $W$ is $Y$ with some 2-handles attached. Dually, $Q$ is obtained from $F$ by attaching tubes along $\Gamma \subset W$, the set of arcs dual to the compressions. We will view $\Gamma$ as a graph in $W$ whose valence one vertices all lie on $F$. At this point $\Gamma$ is just a union of arcs, but as edges of $\Gamma$ are slid onto other edges, higher valence vertices may appear.

We will show that $F$ is incompressible. Assume, in contradiction, that there is a compressing disk $C$ for $F$ in $W$. Isotope $\partial C$ off the disks of $\eta(\Gamma) \cap F$, and use Lemma 1.2 to isotope $\Gamma$ off of $C$. Now $C$ is a compression for $Q$ in $Y$ which is disjoint from a meridian of any tube of $\eta(\Gamma)$, contradicting weak incompressibility of $Q$. Thus $F$ is not compressible into $W$. $F$ is clearly not compressible into $X_-$, since $Q$ has been maximally compressed into $X$, so $F$ is incompressible.

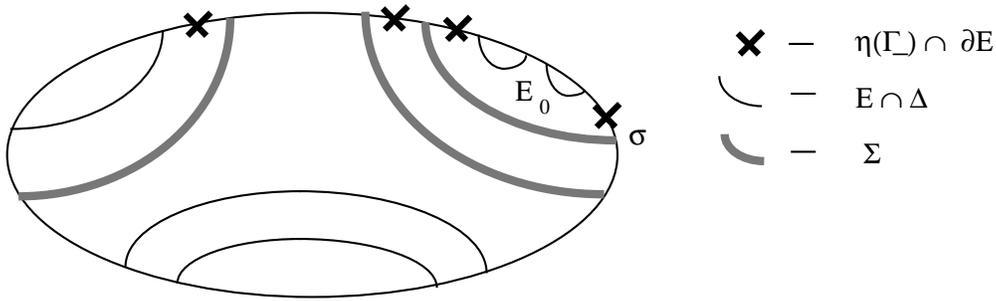

Figure 1.

Isotope $\Delta$ to minimize $|\Delta \cap F|$. Note that if $\Delta \cap F$ contained any closed curves, then each is inessential in $F$ and, by choosing one innermost in $F$ and replacing the disk it bounds in $\Delta$ by a copy of the one it bounds in $F$ (this could be done by an isotopy because $H$ is irreducible so the two disks bound a ball), $|\Delta \cap F|$ could be reduced. Now any compressing disk for $F - \Delta$ in the ball $H - \Delta$ must have a boundary which is inessential in $F$. The disk it bounds in $F$ cannot intersect $\Delta$, since $\Delta \cap F$ contains no closed curves, so the disk in fact lies in $F - \Delta \subset H - \Delta$. We conclude that each component of $F - \Delta$ is incompressible in the ball $H - \Delta$ and so each component of $F - \Delta$ is a disk.

Now $F, \Gamma$, and $\Delta - F$ satisfy the hypothesis of Lemma 1.2, so $\Gamma$ can be isotoped so as to be disjoint from $\Delta$. Regard $\Gamma$ as a disjoint union $\Gamma_- \cup \Lambda$, where $\Gamma_-$ is a graph disjoint from $\Delta$, and $\Lambda$ is a set of pairwise disjoint arcs contained in $\Delta$. We have seen that this can be done with $\Lambda$ empty. The goal is to be able to write $\Gamma$ in this way so that $\Gamma_-$ is empty, i. e. the entire graph $\Gamma$ consists of arcs $\Lambda \subset \Delta$. So let $\Lambda$ and $\Gamma_-$ be chosen to maximize $|\Lambda|$. Then choose a compressing disk $E$ for $Q$ in $Y$ and further slide and isotope $\Gamma_-$ so that the total number of times that $\partial E$ traverses edges of $\Gamma_-$ is minimized. If $\Gamma_-$ is non-empty then $\partial E \cap \eta(\Gamma_-)$ is also non-empty, since, by weak incompressibility of $Q$, $E$ must intersect the meridians of all tubes of $\eta(\Gamma)$. We will show that the assumption that $\Gamma_- \neq \emptyset$ leads inevitably to a contradiction.

Assume first that $E \cap \Delta \neq \emptyset$. Let $\Sigma$ be the set of all arcs of $\Delta \cap E$ such that the two disks into which the arc divides $E$ each intersect $\eta(\Gamma_-)$. If $\Sigma$ is non-empty, let $\sigma$ be outermost in $E$ among arcs of $\Sigma$, and let $E_0$ be the disk of $E - \Delta$ adjacent to $\sigma$ on the $\Sigma$-outermost side. If $\Sigma$ is empty, then all arcs of $\partial E \cap \eta(\Gamma_-)$ are contained in one disk of $E - \Delta$. Call this disk $E_0$ and let $\sigma$ be any component of $\partial E_0 - \partial E$. In either case, the following holds: $\eta(\Gamma_-) \cap \partial E_0 \cap \partial E$ is non-empty, and $\partial E_0 - (\partial E \cup \sigma)$ is a possibly empty collection of arcs each of which cuts off from $E$ a disk disjoint from $\eta(\Gamma_-)$. (See Figure 1.) Let $E'$ be the disk containing $E_0$ which is cut off from $E$ by $\sigma$, so that $E'$ is the union of $E_0$ with a (possibly empty) collection of disks of $E - \Delta$, each of which is disjoint from $\eta(\Gamma_-)$.

Now, if there were an edge $\gamma$ of $\Gamma_-$ with $|\partial E_0 \cap \eta(\gamma)| = 1$, then $E' \supset E_0$ could be used to slide $\gamma$ across until $\gamma$ became $\sigma$. After the slide, the place of $\gamma$ in the decomposition $\Gamma = \Gamma_- \cup \Lambda$ can be changed from $\Gamma_-$ to $\Lambda$, increasing $|\Lambda|$. From this contradiction we conclude that each tube of $\eta(\Gamma_-)$ which intersects $\partial E_0$ intersects it multiple times. Choose a pair of arcs $\gamma_1$ and $\gamma_2$ of $\partial F \cap \eta(\Gamma_-)$ which lie on the same tube $\eta(\gamma)$ of $\eta(\Gamma_-)$ such that no other arc on that same tube lies between them on $\partial E - \sigma$ and so that no other pair coming from a different tube lies between them on $\partial E - \sigma$. Let $\beta$ be the arc in $\partial E' - \sigma$ connecting $\gamma_1$ and $\gamma_2$. (See Figure 2.)

We will now pause for a moment to discuss the case where $E \cap \Delta = \emptyset$. If there were an edge of $\Gamma_-$ whose neighborhood intersected $\partial E$ only once, then $\partial E$ would be a cycle in the ball $H - \Delta$. From this contradiction we can assume that $\partial E$ intersects the neighborhood of each edge of $\Gamma_-$ at least twice. So once again we can choose $\gamma_1, \gamma_2 \in (\partial E \cap \eta(\gamma))$ for some edge $\gamma \in \Gamma_-$, with $\gamma_1$ and $\gamma_2$ connected in $\partial E$ by an arc $\beta$ such that $int(\beta) \cap \eta(\gamma) = \emptyset$ and that there is no edge $\gamma' \in \Gamma_-$ with $|int(\beta) \cap \eta(\gamma')| > 1$. If, for this case, we allow $E_0$ and $E'$ to be alternate names for $E$, the proof can proceed simultaneously with that of the case where $E \cap \Delta \neq \emptyset$:

We claim that the interior of $\beta$ is disjoint from $\eta(\Gamma_-)$. For if it were not, then, by construction, $\beta$ traverses another edge $\gamma' \in \Gamma_-$ exactly once. Let $\beta'$ be an arc in the interior of $E_0$ whose ends coincide with those of $\beta$. Then the subdisk of $E$ containing $\beta$ that $\beta'$ cuts off in $E$ can be used to slide and isotope $\gamma'$ until it coincides with $\beta' \subset E_0$. After the slide, and depending on the relative orientation in $\partial E_0$ of $\gamma_1$ and $\gamma_2$, as identified through $\gamma$, either $\gamma'$ or $\gamma' \cup \gamma$ forms a cycle in the ball $H - \Delta$, contradicting the Corollary 1.5.

We next claim that both ends of $\beta$ lie on the same end $y$ of $\gamma$. For if they lie at opposite ends of $\gamma$ we could construct a cycle in the ball $H - \Delta$ in exactly the same way, after breaking $\gamma$ into two edges, adding a vertex at its midpoint.

It will now be useful to consider the circle $\partial \eta(y) \subset F$. Let $\theta$ and $\theta'$ be the two arcs into which the ends of $\beta$ divide this circle. Now, $H - (F \cup \Delta)$ is obtained from the ball $H - \Delta$ by cutting on the proper disks $F - \Delta$, so $H - (F \cup \Delta)$ is a collection of balls, as is then $H - (F \cup \Delta \cup \eta(\Lambda))$. The closed curve $\beta' \cup \theta$ is parallel (via a subdisk of $E_0$) to a curve on the surface of one such ball and thus $\beta' \cup \theta$ bounds a disk $D$ whose interior is disjoint from $F \cup \Delta \cup \eta(\Lambda)$ (but which may intersect $\Gamma_-$).

Let $K$ be the subdisk of $E$ that bounded by $\beta'$ and $\beta$. Let $F_+$ be the surface obtained from $F$ by attaching tubes along the arcs of $\Lambda$.

**Case 1:** The simple closed curve $\beta \cup \theta \subset F_+$ is essential in $F_+$. Then the surface $F_+$, the graph $\Gamma_-$, and the disk $K_+ = D \cup_{\beta'} K$ satisfy the hypothesis of Lemma 1.2, so, after some slides and isotopies, $\Gamma_-$ can be made disjoint from $K_+$. Since $\partial K_+ = \beta \cup \theta$ is essential in $F_+$ and $F_+$ is obtained from $Q$ by compression, $\partial K_+$ is essential in $Q$. Thus $K_+$ is a compression disk for $Q$ in $Y$ which is disjoint from a meridian disk of any tube of $\eta(\Gamma_-)$. This contradicts weak incompressibility.

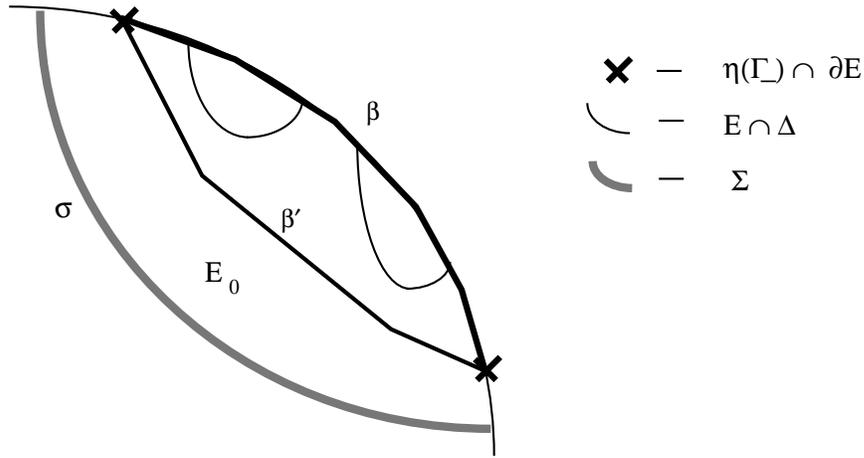

Figure 2.

**Case 2:** $\beta \cup \theta$ is inessential, i. e. bounds a disk $J$ in $F_+$. (See Figure 3.) If $y$ is in $J$ then we can switch from $\theta$ to $\theta'$. So assume that $y$ is not in $J$. If $J$ contains no endpoints of $\Gamma_-$, then the arc $\partial E \cap K$ can be pulled across $J$ to become $\theta$, and then the arc of $\partial E$ which has now become $\gamma_1 \cup \theta \cup \gamma_2$ can be pulled across $\partial \eta(\gamma)$ to eliminate $\gamma_1$ and $\gamma_2$ from $\partial E \cap \eta(\Gamma_-)$. This contradicts our assumption that $\partial E$ minimally intersects the meridia of $\Gamma_-$.

So assume that $J$ contains some endpoints of $\Gamma_-$. Then $D \cup J \cup K$ is a sphere which can be pushed off $F_+$. It bounds a ball $B$ whose interior is disjoint from $F_+$ since $H$ is irreducible and $F$ is incompressible. Choose a small (e. g. disjoint from $\Delta$) collar $D \times I$ of $D$ in $B$, restricting to a small collar of $\theta$ in $J$. Use $B$ to sweep all of $B \cap (\Gamma_- \cup E)$ into $D \times I$ (using most of $J$ to sweep the ends of $\Gamma_-$ into the small collar of $\theta$). Since $D \times I$ is disjoint from $\Delta$, after this move $\Gamma_-$ and $\Delta$ are still disjoint. Now pull the ends of $\Gamma_-$ lying in $D \times I \cap J$ over $\eta(\gamma)$. This move happens entirely in $H - \Delta$, so $\Gamma_-$ and $\Delta$ remain disjoint. We are now back to the case where $J$ is disjoint from $\Gamma_-$, and we can reduce the number of times that $\partial E$ intersects the meridians of the (newly slid) $\Gamma_-$ just as before. Thus we again entcounter a contradiction.

The conclusion is that $\Gamma_- = \emptyset$, as required. □

**Problem:** It is shown in [Sc1] that if $H$ is a solid torus, then $\Lambda$ consists of at most one arc. It is a natural question whether there is in general a bound on $|\Lambda|$ that depends only on the genus of $H$.

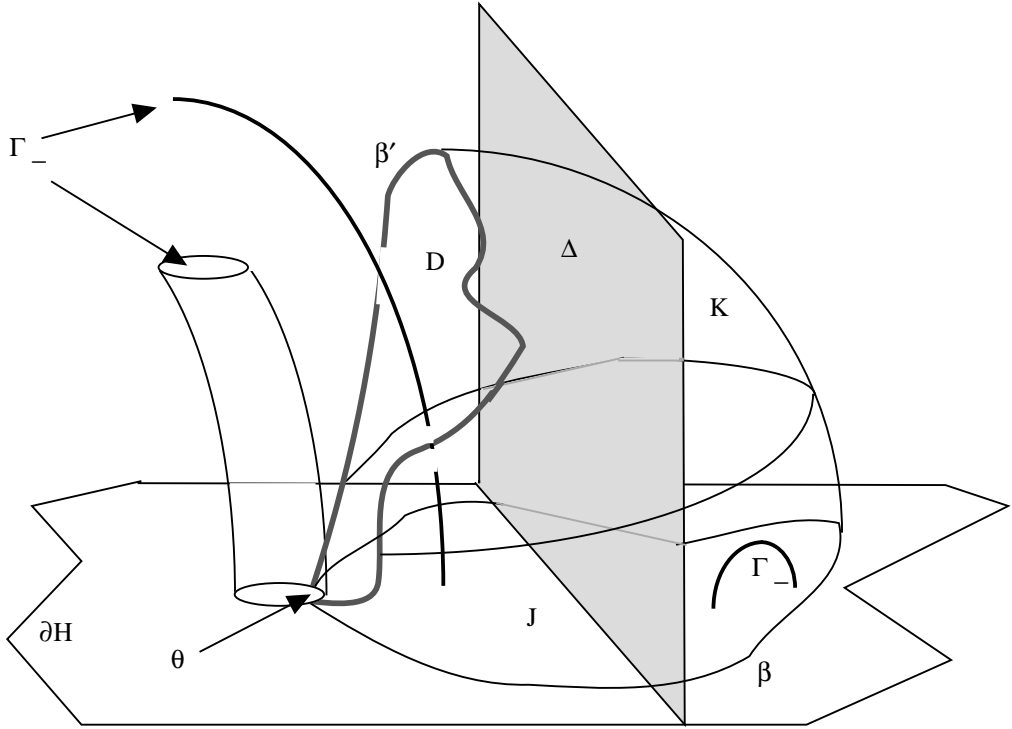

Figure 3.

## 3. How splitting surfaces intersect handlebodies

We intend to use the results of the previous section to understand how a strongly irreducible Heegaard splitting surface $Q$ for a 3-manifold $M$ can intersect a handlebody $H \subset M$.

**Definition 3.1.** *Suppose $(P, \partial P) \subset (M, \partial M)$ is a properly imbedded surface in $M$ and $\alpha$ is an arc on $\partial M$ that intersects $\partial P$ in precisely an end point $x$ of $\alpha$. Then piping $P$ along $\alpha$ is the ambient isotopy of $P$ in $M$ obtained by dragging $x \in P$ along $\alpha$, moving only a small neighborhood of $x$ in $P$.*

Let $X \cup_Q Y$ be a strongly irreducible Heegaard splitting of a 3-manifold $M$ and let $H$ be a handlebody in $M$ with complete collection of meridians $\Delta$.

**Proposition 3.2.** *Suppose $Q \cap \partial H$ consists of curves which are essential in both $\partial H$ and $Q$, the surface $Q \cap H$ compresses into $X \cap H$ and $\partial H \cap Y$ is incompressible in $Y$. Let $F$ be the closed surface obtained from $Q$ by maximally compressing $Q \cap H$ into $X \cap H$, so, dually, $Q$ is obtained by 1-surgery on a graph $\Gamma \subset (H - F)$. Then $\Gamma$ may be slid and isotoped in $M$, $\Delta$ piped in $H$ so that afterwards $F \cap \Delta$ contains no closed curve and $\Gamma$ consists of arcs entirely in $\Delta$ or in $\partial H - \Delta$. Moreover, if $Q \cap H$ also compresses into $Y \cap H$, the sliding and isotopy of $\Gamma$ takes place entirely in $H$.*

**Proof:** : Since $X \cup_Q Y$ is strongly irreducible and $\partial H \cap Y$ is incompressible in $H \cap Y$, it follows that $Q \cap H$ is weakly incompressible. Indeed, an essential closed curve in $Q \cap H$ cannot bound a disk $D_Q$ in $Q$ unless $D_Q$ intersects $\partial H$. Choose such a $D_Q$ to minimize $|D_Q \cap \partial H|$. Consider an innermost (in $D_Q$) such curve of intersection. It must be essential in $\partial H$, so nearby is a compressing disk for $\partial H \cap Y$, contradicting the hypothesis.

The surface $F$ splits $M$ into two 3-manifolds $X^- \subset X$ and $Y^+ \supset Y$. $F \cap H$ is incompressible in $H$. Indeed, a compression into $X_- \cap H$ is impossible by construction. Since $X \cup_Q Y$ is strongly irreducible, the same argument used in the proof of Theorem 2.1 could be used to push $\Gamma$ off of a compressing disk for $F$ in $Y^+$ (the isotopies possibly pushing parts of $\Gamma$ out of $H$). The compressing disk for $F$ then becomes a a compressing disk for $Q$ lying in $Y$, contradicting strong irreducibility. Since $F \cap H$ is incompressible in $H$ it follows as in Theorem 2.1 that $\Delta$ can be isotoped so that $\Delta \cap F$ consists only of arcs and each component of $F - \Delta$ is a disk.

The proof is by induction on the pair $(-\chi(Y \cap \partial H), |\partial H \cap Q|)$, lexicographically ordered. $Y \cap \partial H$ contains no disk components, since curves of intersection are assumed to be essential in $\partial H$. Hence $-\chi(Y \cap \partial H) \geq 0$. Let $E$ be a meridian disk for $Y$, chosen so that $|E \cap \partial H|$ is minimized. If $E$ is disjoint from $\partial H$ then it cannot lie outside $H$, by strong irreducibility, so it must lie inside $H$. In that case the result follows immediately from Theorem 2.1.

If any component of $E \cap \partial H$ is a simple closed curve, then the hypothesis guarantees that an innermost such curve is inessential in $\partial H \cap Y$ and indeed bounds there a disk disjoint from $Q$. An innermost such disk in $\partial H$, could be swapped with the disk its boundary bounds in $E$ to reduce $E \cap \partial H$. So we may as well assume that each component of $E \cap \partial H$ is an arc.

Each arc component of $\partial E \cap (Q - \partial H)$ is essential in the bounded surface $Q - \partial H$. For if not, then an outermost (in $Q - \partial H$) arc of intersection could be slid to the other side of $\partial H$, either connecting two arc components of $E \cap \partial H$ or creating a single closed component which can then be eliminated as above.

Consider an $E$-outermost arc $\delta$ of $E \cap \partial H$, cutting off of a disk $E'$ from $E$. If $E' \subset M - H$ then use $E'$ to isotope an essential arc of $Q - H$ into $Q$. This move cuts $Y \cap \partial H$ along an arc, and so decreases $-\chi(Y \cap \partial H)$ by one. Any closed curve in $\partial H \cap Y$ that is essential after the cut was essential before, so $\partial H \cap Y$ remains incompressible in $Y$. If the $\partial$-compression creates an inessential curve in $\partial H$ (necessarily bounding a disk in $\partial H \cap Y$) then that curve must also be inessential in $Q$, by strong irreducibility, so can be isotoped into $H$. The net result is to isotope an annulus component $A$ of $Q - H$ across a parallel annulus component of $Y \cap \partial H$ and into $H$. In any case, this reduces the pair $(-\chi(Y \cap \partial H), |\partial H \cap Q|)$. Moreover $F \cap H$ remains incompressible in $X^- \cap H$. To see that this is true even in the case when an annulus $A$ is pushed across, note that the isotopy of $A$ can be undone by pushing the core of $A$ to $\partial H$ and across along an annulus $B \subset X^-$ that has one end on $\partial H$ and the other on $Q \cap H$. If $F \cap H$ compressed in $X^- \cap H$ then the compressing disk would necessarily intersect $B$ and all components of intersection would be arcs in $B$ with both ends on $F \cap H$. An outermost (in $B$) such arc could be used to alter the compressing disk, eventually creating a compressing disk for $F \cap H$ in $H$ that is disjoint from $B$, a contradiction. So ultimately, all the inductive hypotheses still apply, $\Gamma$ may now be slid and isotoped into place, either in a collection of meridian disks $\Delta$ or in $\partial H - \Delta$. Then push the annulus $A$ back across $\partial H$. Since the push is across $X^- \cap \partial H$ this move does not affect the positioning of arcs of $\Gamma$, which lie in $Y^+$.

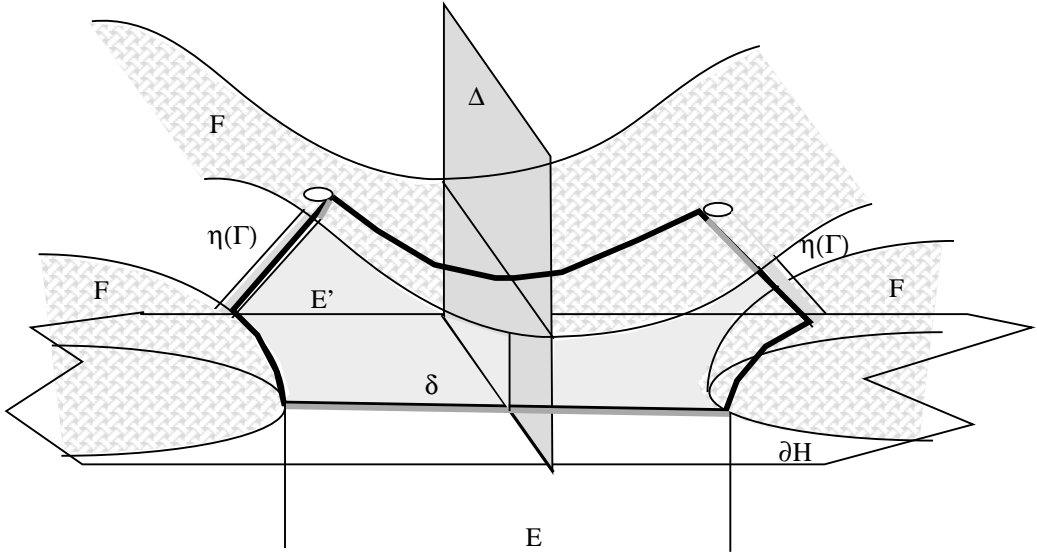

Figure 4.

The hard case is when $E' \subset H$, since then $\partial E'$ may run along parts of $\Gamma$. (See Figure 4.) The idea is to mimic the proof of Theorem 2.1, using $E'$ instead of $E$. There are two important differences: The subarc $\delta \subset \partial E'$ lies on $\partial H$ (so components of $E' \cap \Delta$ may have boundary points on $\delta \subset \partial H$). And there is no guarantee that $\partial E'$ passes over all (or indeed any) of the arcs of $\Gamma$.

First mimic the proof of Lemma 1.2 to slide and isotope edges of $\Gamma$ through $M$ (perhaps temporarily leaving $H$) in order to minimize transverse intersections of $\Gamma \cap \Delta$. The situation is complicated by the fact that arcs of $E \cap \Delta$ may have ends on $\delta$. Just as in the proof of Lemma 1.2, there is a point $w \in \Delta \cap \Gamma$ so that any arc of intersection of $E$ with $\Delta$ that is incident to $w$ (and there must be at least two, by strong irreducibility) has its other end either on $\partial H$ or on $F$ or on another point of $\Delta \cap \Gamma$.

Consider an arc of $\Delta \cap E$ which has an end at $w$. The other end either also lies on $\partial E$ or on an arc component $\kappa$ of $\partial H \cap E$. In either case, it, and possibly a subarc of $\kappa$, cut off a sector from $E$ that is disjoint from $\kappa$. Choose an outermost such sector in $E$. This is a disk whose boundary consists of a subarc of $\partial E$ whose interior is disjoint from $w$ (since the sector is outermost of those that intersect $w$), an arc component $\epsilon$ of $E \cap \Delta$, and possibly a subarc $\kappa_0$ of a component of $\partial H \cap E$. If both ends of $\epsilon$ lie on $\partial E$, so $\kappa_0$ is not used, then the sector can be used to slide and isotope a segment of the edge containing $w$ so that it lies on $\Delta$ (the first case) and so can be isotoped off of $\Delta$, reducing $\Delta \cap \Gamma$ by at least one. In the other case, first use $\kappa_0$ to sequentially pipe the disks $\Delta$ across the end of $\kappa_0$ at $\partial E$. After this maneuver, the ends of $\epsilon$ lie on $\partial E$ and we are done as before.

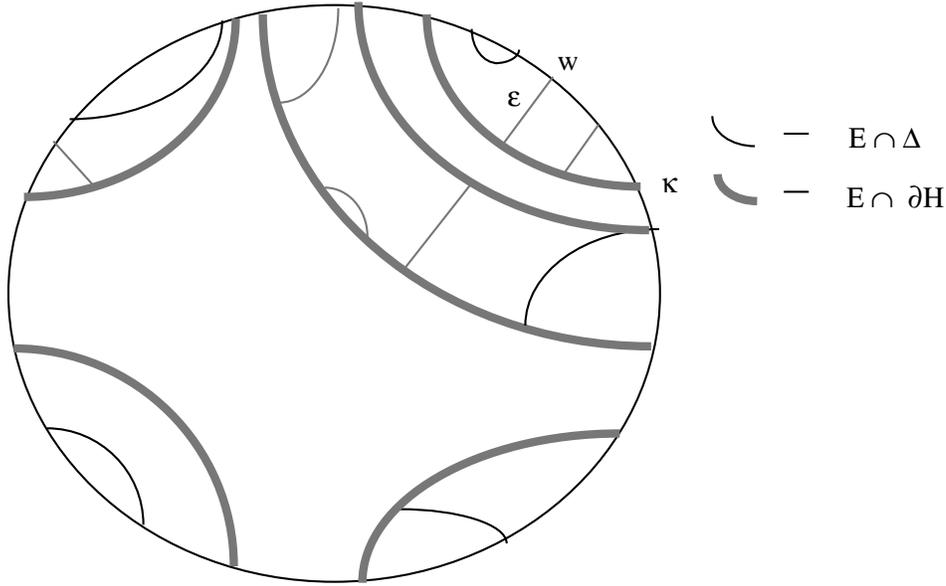

Figure 5.

Now proceed much as in the proof of Theorem 2.1: Consider $\Gamma$ as the union of a graph $\Gamma_-$ disjoint from $\Delta$ and a collection of arcs $\Lambda$ which lie entirely in $\Delta$ or entirely in $\partial H - \Delta$. We start with $\Lambda = \emptyset$ and the goal is to isotope and slide to achieve $\Gamma_- = \emptyset$. First maximize $|\Lambda|$. Let $E' \subset H$ be an outermost sector cut off by $\partial H$ from $E$ as above.

**Case 1:** $\partial E'$ traverses (i. e. intersects the meridian) of every edge in $\Gamma_-$.

In this case, apply a variant of the argument in Theorem 2.1, with minor changes made at the beginning of the argument: After piping the disks $\Delta$ along subarcs of the arc $\delta = \partial E' \cap \partial H$, we can assume that $E' \cap \Delta$ consists entirely of arcs with both ends on $Q \cap H$.

Let $\Sigma$ be those arcs in $\Delta \cap E'$ which have the property that each cuts off a subsector which doesn't contain $\delta$, but does intersect $\Gamma_-$. If $\Sigma$ is non-empty, let $\sigma$ be outermost in $E'$ among arcs of $\Sigma$, and let $E_0$ be the disk of $E' - \Delta$ adjacent to $\sigma$ on the $\Sigma$-outermost side. If $\Sigma$ is empty, then all arcs of $\partial E' \cap \eta(\Gamma_-)$ are contained in the component of $E' - \Delta$ that contains $\delta$. Call this disk $E_0$ and let $\sigma = \delta$. In either case, the following holds: $\eta(\Gamma_-) \cap \partial E_0 \cap \partial E$ is non-empty, and $\partial E_0 - (\partial E' \cup \sigma)$ is a possibly empty collection of arcs each of which cuts off from $E'$ a disk disjoint from $\eta(\Gamma_-)$. Let $E''$ be the disk containing $E_0$ which is cut off from $E$ by $\sigma$, so that $E''$ is the union of $E_0$ with a (possibly empty) collection of disks of $E' - \Delta$, each of which is disjoint from $\eta(\Gamma_-)$.

Now, if there were an edge $\gamma$ of $\Gamma_-$ with $|\partial E_0 \cap \eta(\gamma)| = 1$, then $E'' \supset E_0$ could be used to slide $\gamma$ across until $\gamma$ became $\sigma$. After the slide, the place of $\gamma$ in the decomposition $\Gamma = \Gamma_- \cup \Lambda$ can be changed from $\Gamma_-$ to $\Lambda$, increasing $|\Lambda|$. From this contradiction we conclude that each tube of $\eta(\Gamma_-)$ which intersects $\partial E_0$ intersects it multiple times. Now the proof in this case concludes exactly as in Theorem 2.1.

**Case 2:** $\partial E'$ fails to traverse some edge in $\Gamma_-$.

Let $\Gamma'$ be the subgraph of $\Gamma_-$ whose edges are disjoint from $\partial E'$. First $\partial$-compress $Q$ to $\partial H$ via $E'$. This decreases $-\chi(Y \cap \partial H)$ but leaves $\Gamma'$ untouched. The inductive hypothesis then allows $\Gamma'$ to be slid into place. We need to ascertain that after this process, one can undo the $\partial$-compression in a way that returns $\Lambda$ to its former position, so that $\Lambda$ is, in the end, augmented by $\Gamma'$, completing the inductive step. In order to accomplish this convincingly, we will add a further level of induction, on the number of arcs in $E' \cap \Delta$.

If $E' \cap \Delta = \emptyset$ then first $\partial$-compress $Q \cap H$ along $E'$, decreasing $-\chi(Y \cap \partial H)$. Then apply the inductive step to isotope $\Gamma'$ until it consists of arcs in $\Delta \cap Y^+$ and $\partial H \cap Y^+$. The $\partial$-compression via $E'$ can then be undone via a $\partial$-compressing disk that lies in $X^- \cap (M - H)$. Since the $\partial$-compression only crosses $X^- \cap \partial H$ it has no effect on $\Gamma' \cap \partial H \subset Y^+$.

If $E' \cap \Delta \neq \emptyset$, begin with two technical simplifications. Pipe $\Delta$ along end segments of the arc $\delta$ so that afterwards, each arc of $E' \cap \Delta$ has both ends on $F$. Then an outermost sector of $E'$ cut off by $\Delta$ can be used to isotope $Q \cap H$ to reduce $|E' \cap \Delta|$. After this move, the inductive assumption allows us to isotope $\Gamma'$ until it lies in $\Delta \cap Y^+$ or $(\partial H - \Delta) \cap Y^+$. Then the first move of $Q \cap H$ across $Y \cap \Delta$ can be undone by a push across $X^- \cap \Delta$, and this has no effect on the arcs of $\Gamma'$ that have been moved to $Y^+ \cap \Delta$ or $(\partial H - \Delta) \cap Y^+$. $\square$

## 4. Intersections as unknotted surfaces

**Definition 4.1.** *A* spine *of a surface $P$ is a 1-complex $\Sigma$ properly imbedded in $P$ so that each component of $P - (\Sigma \cup \partial P)$ is an open disk.*

**Definition 4.2.** *A surface $P, \partial P$ properly imbedded in a handlebody $H, \partial H$ is* unknotted *if for some spine $\Sigma$ of $P$, $H - \Sigma$ is also a handlebody.*

The notion of unknotted surface in a compression body was introduced (without the terminology) in [RS2, Section 4]. It is shown there that in fact the choice of spine is irrelevant; if $P$ is unknotted using one spine, it is unknotted for all spines. Furthermore it is shown that incompressible surfaces in $H$ are always unknotted, as are weakly incompressible bicompresible surfaces without closed component. Indeed, the entire discussion there is generalized to surfaces in compression bodies.

The terminology allows us to incorporate Theorem 2.1 and Proposition 3.2 into a general statement about how a splitting surface $Q$ of a strongly irreducible splitting $X \cup_Q Y$ intersects a handlebody $H \subset M$.

**Theorem 4.3.** *Let $X \cup_Q Y$ be a strongly irreducible Heegaard splitting of a 3-manifold $M$ containing a handlebody $H$. Suppose each curve in $Q \cap \partial H$ is essential in both $Q$ and $\partial H$ and the surfaces $\partial H \cap X$ and $\partial H \cap Y$ are incompressible in $X$ and $Y$ respectively. Then $Q \cap H$ is unknotted in $H$.*

**Proof:** If $Q \cap H$ is incompressible in $H$, this follows immediately from [RS2, Prop. 4.2]. If $Q \cap H$ compresses into $X \cap H$, say, the result follows from 3.2, as follows.

Choose a complete meridian system $\Delta$ of $H$ with the property that each ball in $H - \Delta$ is incident to each disk in $\Delta$ on at most one side. An early step in the proof of Proposition 3.2 is to isotope so that each component of $F - \Delta$ is a disk (so in particular $F \cap \Delta$ is a spine of $F$), and this property is not affected by the piping of $\Delta$. The addition to $F$ of those tubes of $\Lambda$ that lie in $\Delta$ has the effect of banding together these disks in $H - \Delta$, and, by [Fr, Lemma 1.1] (the precursor to Lemma 1.5) this banding of disks cannot create a more complicated surface than a disk. So even after these tubes are attached the intersection $\Sigma$ of the surface with $\Delta$ is still a spine of the surface. The remaining tubes that need to be attached to obtain $Q \cap H$ are tubes parallel to arcs in $\partial H - \Delta$. For each such tube, augment $\Sigma$ by an arc running once through the tube and down to $\partial H$ at each end. The result in the complement of $\Sigma$ is again to band two disks together to get another disk. At the end of the process, $\Sigma$ remains a spine of the surface, which is now all of $Q \cap H$. Moreover, each arc of $\Sigma$ is parallel to an arc in $\partial H$, either because it lies in $\Delta$ or by construction. Hence $H - \eta(\Sigma)$ is a handlebody. □

MATT JONES, DEPARTMENT OF COGNITION AND PERCEPTION, UNIVERSITY OF MICHIGAN, ANN ARBOR, MI 48109

MARTIN SCHARLEMANN, MATHEMATICS DEPARTMENT, UNIVERSITY OF CALIFORNIA, SANTA BARBARA, CA USA
*E-mail address*: `mgscharl@math.ucsb.edu`